\documentclass[12pt]{article}
\usepackage{amsfonts,amsmath,amsthm, amssymb}
\usepackage{latexsym, euscript, epic, eepic}
\renewcommand{\theequation}{\thesection.\arabic{equation}}
\newtheorem{cro}{Corollary}[section]

\newtheorem{thm}{Theorem}[section]
\newtheorem{lem}{Lemma}[section]
\newtheorem{rem}{\bf Remark}[section]

\begin{document}
\title{Packing dimensions of the divergence points of self-similar
measures with the open set condition
%\footnotetext {* Corresponding author}
 }
\author{Xiaoyao Zhou $^\dag,$  Ercai Chen$^{\dag \ddag}$\\
\small \it $\dag$ School of Mathematical Science, Nanjing Normal University,\\
\small \it Nanjing 210097, Jiangsu, P.R.China.\\
\small \it $\ddag$ Center of Nonlinear Science, Nanjing University\\
\small \it Nanjing 210093, Jiangsu, P.R.China.\\
\small \it e-mail:
\small \it\noindent zhouxiaoyaodeyouxian@126.com\\
\small \it \noindent ecchen@njnu.edu.cn}
\date{}
\maketitle

\begin{center}
\begin{minipage}{120mm}
{\small {\bf Abstract.} Let $\mu$ be the self-similar measure
supported on the self-similar set $K$ with open set condition. In
this article, we discuss the packing dimension of the set $\{x\in K:
A(\frac{\log\mu(B(x,r))}{\log r})=I\}$ for $I\subseteq\mathbb{R},$
where $A(\frac{\log\mu(B(x,r))}{\log r})$ denotes the set of
accumulation points of $\frac{\log\mu(B(x,r))}{\log r}$ as
$r\searrow0.$ Our main result solves the conjecture about packing
dimension posed by Olsen and Winter \cite{OlsWin} and generalizes
the result in \cite{BaeOlsSni}.   }
\end{minipage}
\end{center}

\vskip0.5cm
%%%%%%%%%%%%%%%%%%%%%%%%%%%%%%%%%%%%%%%%%%%%%%%%%%%%%%%%%%%%%%%%%%%%%%%%%%%%%%%%%%%%%%%%%%%%%%%%%%%%%%%%%%%%%%%%%%%
%%%%%%%%%%%%%%%%%%%%%%%%%%%%%%%%%%%%%%%%%%%%%%%%%%%%%%%%%%%%%%%%%%%%%%%%%%%%%%%%%%%%%%%%%%%%%%%%%%%%%%%%%%%%%%%%%%%
%%%%%%%%%%%%%%%%%%%%%%%%%%%%%%%%%%%%%%%%%%%%%%%%%%%%%%%%%%%%%%%%%%%%%%%%%%%%%%%%%%%%%%%%%%%%%%%%%%%%%%%%%%%%%%%%%%%
%%%%%%%%%%%%%%%%%%%%%%%%%%%%%%%%%%%%%%%%%%%%%%%%%%%%%%%%%%%%%%%%%%%%%%%%%%%%%%%%%%%%%%%%%%%%%%%%%%%%%%%%%%%%%%%%%%%
%%%%%%%%%%%%%%%%%%%%%%%%%%%%%%%%%%%%%%%%%%%%%%%%%%%%%%%%%%%%%%%%%%%%%%%%%%%%%%%%%%%%%%%%%%%%%%%%%%%%%%%%%%%%%%%%%%%
\section{Introduction and statement of results}
Let $S_i:\mathbb{R}^d\to\mathbb{R}^d$ for $i=1,\cdots,N$ be
contracting similarities with contraction ratios $r_i\in (0,1)$ and
let $(p_1,\cdots, p_N)$ be a probability vector. Let $K$ and $\mu$
be the self-similar set and the self-similar measure associated with
the list $(S_1,\cdots, S_N, p_1,\cdots, p_N),$ namely, $K$ is the
unique non-empty compact subset of $\mathbb{R}^d$ such that

\begin{equation*}
K=\bigcup\limits_i S_i(K)
\end{equation*}
and $\mu$ is the unique Borel probability measure on $\mathbb{R}^d$
such that

\begin{equation}\label{formula1.1}
\mu=\sum\limits_i p_i\mu\circ S_i^{-1}.
\end{equation}
As we all known that supp$\mu=K.$ The Open Set Condition (OSC) is
fulfilled means that if there exists an open non-empty and bounded
subset $U$ of $\mathbb{R}^d$ such that $\bigcup\limits_i S_i
U\subseteq U$ and $S_i U\cap S_j U=\emptyset$ for all $i,j$ with
$i\neq j.$

During the past 20 years the multifractal structure of $\mu$ has
attracted considerable attention. Multifractal analysis refers to
the study of the fractal geometry of the set
$\left\{x\in\mathbb{R}^d:\lim\limits_{r\to0}\frac{\log\mu(B(x,r))}{\log
r}=\alpha\right\}.$ Define the Hausdorff multifractal spectrum
$f_H(\alpha)$ of $\mu$ and the packing multifractal spectrum
$f_P(\alpha)$ of $\mu$ as below. For $\alpha\geq0,$ put

\begin{equation}
f_H(\alpha)=\dim_H\left\{x\in\mathbb{R}^d:\lim\limits_{r\to0}\frac{\log\mu(B(x,r))}{\log
r}=\alpha\right\},
\end{equation}
and

\begin{equation}
f_P(\alpha)=\dim_P\left\{x\in\mathbb{R}^d:\lim\limits_{r\to0}\frac{\log\mu(B(x,r))}{\log
r}=\alpha\right\}.
\end{equation}

Define the function $\beta(q):\mathbb{R}\to\mathbb{R}$ by

\begin{equation}\label{formula1.4}
\sum\limits_{i=1}^N p_i^qr_i^{\beta(q)}=1.
\end{equation}
The Legendre transformation of a real valued function
$\varphi:\mathbb{R}\to\mathbb{R}$ is given by
$\varphi^*(x)=\inf\limits_y(xy+\varphi(y)).$
 Arbeiter and Patzschke
\cite{ArbPat}\cite{Pat} succeed in computing the multifractal
spectra $f_H(\alpha)$ and $f_P(\alpha)$ under the OSC.

\begin{thm}{\rm(see \cite{ArbPat} \cite{Pat})}
Assume that the OSC is satisfied. Then the multifractal spectra
$f_H(\alpha)$ and $f_P(\alpha)$ are given by

\begin{equation*}
f_H(\alpha)=f_P(\alpha)=\beta^*(\alpha),{\rm~~for~} \alpha\geq0.
\end{equation*}
\end{thm}

\begin{rem}
Let $s=\dim_HK=\dim_PK.$

{\rm\noindent(1)} $\beta(q)$ is a strong decreasing function and
$\beta(1)=0,\beta(0)=s.$

{\rm\noindent(2)} It is either $\beta(q)=s(1-q)$ for
$p_i=r_i^s,i=1,2\cdots,N$ or $\beta(q)$ is convex.

{\rm\noindent(3)} $\alpha $ can considered as a function about $q,$
and $\alpha=\alpha(q)$ is either $\forall q, \alpha(q)=s$ or a
strong decreasing function.
\end{rem}

However, the limit $\lim\limits_{r\to0}\frac{\log\mu(B(x,r))}{\log
r}$ may not exist. Points $x$ for which the limits do not exist are
called divergence points. Strichartz \cite{Str} showed that

\begin{equation*}
\mu\left\{x\in K:\liminf\limits_{r\to0}\frac{\log\mu(B(x,r))}{\log
r}<\limsup\limits_{r\to0}\frac{\log\mu(B(x,r))}{\log r}\right\}=0.
\end{equation*}
However, Barreira and Schmeling \cite{BarSch} and Chen and Xiong
\cite{CheXio} proved that

\begin{equation*}
\left\{x\in K:\liminf\limits_{r\to0}\frac{\log\mu(B(x,r))}{\log
r}<\limsup\limits_{r\to0}\frac{\log\mu(B(x,r))}{\log r}\right\}
\end{equation*}
has full Hausdorff dimension for $(p_1,\cdots, p_N)\neq
(r_1^s,\cdots,r_N^s),$ where $s=\dim_H K.$ Chen and Xiong obtained
the above result under the strong separation condition. Xiao, Wu and
Gao \cite{XiaWuGao} proved that Chen and Xiong's result remains
valid under the OSC. These results shows that the set of divergence
points has an extremely rich and intricate fractal structure. Hence,
Dividing divergence points into different level sets like
multifractal analysis has attracted  an enormous interest in the
mathematical literature.

\noindent For $x\in K,$ set
\begin{equation*}
A(D(x))=\{y\in(0,+\infty):\lim\limits_{k\to\infty}D_{r_k}(x)=y
{\rm~for ~some~}\{r_k\}_k\searrow0\},
\end{equation*}
where $D_r(x)=\frac{\log\mu(B(x,r))}{\log r}.$

We now state Olsen and Winter's results and I. Baek, L. Olsen and N.
Snigireva's result as follows.

\begin{thm}{\rm\cite{OlsWin}}\label{thm1.2}
Assume that the strong separation condition is satisfied. Let
$\alpha_{\min}=\min\limits_i \frac{\log p_i}{\log r_i}$ and
$\alpha_{\max}=\max\limits_i \frac{\log p_i}{\log r_i}.$

\noindent {\rm (1)} If $I\subseteq \mathbb{R}$ is not a closed
subinterval of $[\alpha_{\min},\alpha_{\max}],$ then

\begin{equation*}
\{x\in K: A(D(x))=I\}=\emptyset.
\end{equation*}

 \noindent{\rm (2)} If $I\subseteq \mathbb{R}$ is a closed
subinterval of $[\alpha_{\min},\alpha_{\max}],$ then

\begin{equation*}
\dim_H\{x\in K: A(D(x))=I\}=\inf\limits_{\alpha\in
I}\beta^*(\alpha).
\end{equation*}
\end{thm}

\begin{thm}{\rm\cite{OlsWin}}\label{thm1.3}
Assume that the strong separation condition is satisfied. If
$C\subset\mathbb{R}$ is an arbitrary subset of $\mathbb{R},$ then

\begin{equation*}
\dim_H\{x\in K: A(D(x))\subseteq I\}=\dim_P\{x\in K:
A(D(x))\subseteq I\}=\sup\limits_{\alpha\in C}\beta^*(\alpha).
\end{equation*}
\end{thm}

\begin{thm}{\rm\cite{BaeOlsSni}}\label{thm1.4}
Assume that the strong separation condition is satisfied.

\noindent {\rm (1)} If $I\subseteq \mathbb{R}$ is not a closed
subinterval of $[\alpha_{\min},\alpha_{\max}],$ then

\begin{equation*}
\{x\in K: A(D(x))=I\}=\emptyset.
\end{equation*}

 \noindent{\rm (2)} If $I\subseteq \mathbb{R}$ is a closed
subinterval of $[\alpha_{\min},\alpha_{\max}],$ then

\begin{equation*}
\dim_P\{x\in K: A(D(x))=I\}=\sup\limits_{\alpha\in
I}\beta^*(\alpha).
\end{equation*}
\end{thm}
Olsen and Winter \cite{OlsWin} conjectured that the results in
Theorem \ref{thm1.2} and Theorem \ref{thm1.3} remain valid even if
the strong separation condition is replaced by the OSC. J. Li, M. Wu
and Y. Xiong \cite{LiWuXio} proved that Theorem \ref{thm1.2} remains
valid even if the strong separation condition is replaced by OSC.
This article solves Olsen and Winter's conjecture about packing
dimension positively. Namely, we will prove that Theorem
\ref{thm1.3} and Theorem \ref{thm1.4} remain valid even if the
strong separation condition is replaced by OSC. More precisely, this
paper has the following main results.

\begin{thm}\label{thm1.5}
Assume that the OSC is satisfied.

\noindent {\rm (1)} If $I\subseteq \mathbb{R}$ is not a closed
subinterval of $[\alpha_{\min},\alpha_{\max}],$ then

\begin{equation*}
\{x\in K: A(D(x))=I\}=\emptyset.
\end{equation*}

 \noindent{\rm (2)} If $I\subseteq \mathbb{R}$ is a closed
subinterval of $[\alpha_{\min},\alpha_{\max}],$ then

\begin{equation*}
\dim_P\{x\in K: A(D(x))=I\}=\dim_P\{x\in K: A(D(x))\subseteq
I\}=\sup\limits_{\alpha\in I}\beta^*(\alpha).
\end{equation*}
\end{thm}

\begin{cro}\label{corollary}
Assume that the OSC is satisfied. If $I\subseteq \mathbb{R}$ is a closed
subinterval of $[\alpha_{\min},\alpha_{\max}],$ then

\begin{equation*}
\dim_H\{x\in K: A(D(x))\subseteq I\}=\sup\limits_{\alpha\in
I}\beta^*(\alpha).
\end{equation*}
\end{cro}

\begin{rem}
Here, singleton is regarded as a special closed subinterval.
\end{rem}
%%%%%%%%%%%%%%%%%%%%%%%%%%%%%%%%%%%%%%%%%%%%%%%%%%%%%%%%%%%%%%%%%%%%%%%%%%%%%%%%%%%%%%%%%%%%%%%%%%%%%%%%%%%%%%%%%%%
%%%%%%%%%%%%%%%%%%%%%%%%%%%%%%%%%%%%%%%%%%%%%%%%%%%%%%%%%%%%%%%%%%%%%%%%%%%%%%%%%%%%%%%%%%%%%%%%%%%%%%%%%%%%%%%%%%%
%%%%%%%%%%%%%%%%%%%%%%%%%%%%%%%%%%%%%%%%%%%%%%%%%%%%%%%%%%%%%%%%%%%%%%%%%%%%%%%%%%%%%%%%%%%%%%%%%%%%%%%%%%%%%%%%%%%
%%%%%%%%%%%%%%%%%%%%%%%%%%%%%%%%%%%%%%%%%%%%%%%%%%%%%%%%%%%%%%%%%%%%%%%%%%%%%%%%%%%%%%%%%%%%%%%%%%%%%%%%%%%%%%%%%%%
%%%%%%%%%%%%%%%%%%%%%%%%%%%%%%%%%%%%%%%%%%%%%%%%%%%%%%%%%%%%%%%%%%%%%%%%%%%%%%%%%%%%%%%%%%%%%%%%%%%%%%%%%%%%%%%%%%%
\section{Preliminaries}
We begin by introducing the definition of packing dimension, which
is referred to \cite{Fal1}. Let $0\leq s<\infty. $ For
$A\subseteq\mathbb{R}^d$ and $0<\delta<\infty,$ put
$P^s_\delta(A)=\sup\sum\limits_{i}d(B_i)^s,$ where $d(B_i)$ denotes
the diameter of $B_i$ and the supremum is taken over all disjoint
families of closed balls $\{B_1,B_2,\cdots\}$ such that $d(B_i)\leq
\delta$ and the centres of the $B_i'$s are in $A.$ Define

\begin{equation*}\begin{split}
&P^s(A)=\lim\limits_{\delta\to0}P^s_\delta(A),\\&
\mathcal{P}^s(A)=\inf\left\{\sum\limits_{i=1}^\infty P^s(A_i) :
A=\bigcup\limits_{i=1}^ \infty A_i\right\},
\end{split}\end{equation*} and

\begin{equation*}\begin{split}
\dim_P
A&=\inf\{s:\mathcal{P}^s(A)=0\}=\inf\{s:\mathcal{P}^s(A)<\infty\}\\&
=\sup\{s:\mathcal{P}^s(A)>0\}=\sup\{s:\mathcal{P}^s(A)=\infty\}.
\end{split}\end{equation*}

Fix $N\in \mathbb{N},$ for $n=1,2,\cdots,$ set $\Sigma^n=\{1,\cdots,
N\}^n$  that is, $\Sigma^n$ is the family of all lists $u=u_1\cdots
u_n$ of length $n$ with entries $u_j\in\{1,\cdots, N\}$ and let
$\Sigma^*=\bigcup\limits_n\Sigma^n.$ For $u\in\Sigma^n,$ write
$|u|=n$ for the length of $u.$ For $u\in \Sigma^*, u^-$ is the word
obtained from $u$ by dropping the last letter. Finally, if
$u=u_1\cdots u_n\in \Sigma^n,$ write $S_u=S_{u_1}\circ \cdots \circ
S_{u_n}, p_u=p_{u_1}\cdots p_{u_n}$ and $r_u=r_{u_1}\cdots r_{u_n},$
put $K_u=S_uK$ and set $r_{\min}=\min\limits_{1\leq i\leq N}r_i.$

Let $\mu$ be a Borel probability measure on $\mathbb{R}^d$ admitting
compact support. For any open set $V\subseteq \mathbb{R}^d$ with
$\mu(V)>0,$ put

\begin{equation*}
\Theta_V(q;
r)=\sup\sum\limits_i\mu(B(x_i,r))^q,~~r>0,q\in\mathbb{R},
\end{equation*}
where the supremum is taken over all families of disjoint closed
balls $\{B(x_i,r)\}_i$ contained in $V$ with $x_i\in {\rm supp}
\mu.$ The $L^q$-spectrum $\tau_V(q)$ of $\mu$ on $V$ is defined by

\begin{equation}
\tau_V(q)=\liminf_{r\to0}\frac{\log \Theta_V(q;r)}{-\log r}.
\end{equation}
Particularly, for $V=\mathbb{R}^d,$ simplify
$\Theta_{\mathbb{R}^d}(q;r)$ as $\Theta(q;r),$
$\tau_{\mathbb{R}^d}(q)$ as $\tau(q).$ Peres and Solomyak
\cite{PerSol} showed that for any self-similar measure $\mu,$ the
limit $\lim\limits_{r\to0}\frac{\log\Theta(q;r)}{\log r}$ exists for
$q\geq0.$ In particular, Lau \cite{Lau} obtained that if the OSC is
satisfied, then $\tau(q)=\beta(q)$ for any $q>0,$ where $\beta(q)$
is defined as in (\ref{formula1.4}). Riedi \cite{Rie} proved that
(\ref{formula1.4}) remains valid for $q\in\mathbb{R}.$

If $(S_i)_{i=1}^N$ satisfies the OSC, then there exists an open,
bounded and non-empty set $U$ such that
$\bigcup\limits_{i}S_i(U)\subset U, U\cap K\neq \emptyset$ and
$S_i(U)\cap S_j(U)=\emptyset$ for all $i\neq j.$ It is simple to
check that $S_u K\cap S_{u'}U=\emptyset,$ for any $u\neq
u'\in\Sigma^n.$ Therefore, one can find an open ball
$U_0=U(x_0,r_0)\subset U$ with $x_0\in K,$ where $U(x_0,r_0)$ is the
open ball of radius $r_0$ centered at $x_0.$ And we fix $U_0$ in
this article.

The following two lemmas can be found in \cite{LiWuXio}.
\begin{lem}{\rm\cite{LiWuXio}}\label{lem2.1}
If $\mu$ is a self-similar measure supported on $K$ and satisfies
the OSC, then $\tau_{U_0}(q)=\tau(q)=\beta(q)$ for any
$q\in\mathbb{R}.$
\end{lem}

\begin{lem}{\rm\cite{FenLau}\cite{LiWuXio}}\label{lem2.2}
If we choose any $q\in \mathbb{R}$ and set $\alpha=-\beta'(q),$ then
for any $\delta,\eta>0,$ there exist $d\in (0,\eta), k\geq
d^{-\beta^*(\alpha)+\delta(|q|+1)}$ and $u_1,\cdots, u_k\in
\Sigma^*$ such that

\noindent{\rm (a)} $d^{1+\delta}\leq r_{u_i}\leq d^{1-\delta} $ for
all $1\leq i\leq k.$

\noindent{\rm (b)} $S_{u_i}(U(x_0, 4r_0))$ are disjoint subsets of
$U_0.$

\noindent{\rm (c)} $d^{\alpha+3\delta}\leq p_{u_i}\leq
d^{\alpha-3\delta}$ for all $1\leq i\leq k.$
\end{lem}
To get the lower bound estimate of the packing dimension in Theorem
\ref{thm1.5}, we need to construct a Moran set. Let's present the
definition of the Moran set and some results of it. Fix a closed
ball $B\subset \mathbb{R}^d.$ Let $\{n_k\}_{k\geq1}$ be a sequence
of positive integers. Let $D=\bigcup\limits_{k\geq0}D_k$ with
$D_0=\{\emptyset\}$ and $D_k=\{\omega=(j_1 j_2\cdots j_k):1\leq
j_i\leq n_i, 1\leq i\leq k\}.$ Suppose that
$\mathcal{G}=\{B_\omega:\omega\in D\}$ is a collection of closed
balls of radius $r_\omega$ in $\mathbb{R}^d.$ We say that
$\mathcal{G}$ fulfills the Moran structure provided it satisfies the
following conditions:

\begin{itemize}
  \item $B_\emptyset=B, B_{\omega j}\subset B_{\omega}$ for any
$\omega\in D_{k-1}, 1\leq j\leq n_k;$
  \item $B_\omega\cap B_{\omega'}=\emptyset$ for $\omega,
\omega'\in D_k$ with $\omega\neq\omega';$
  \item $\lim\limits_{k\to\infty}\max\limits_{\omega\in
D_k}r_\omega=0;$
  \item For all $\omega\eta\neq \omega'\eta, \omega,\omega'\in
D_m, \omega\eta, \omega\eta\in D_n, m< n,$

\begin{equation*}
\frac{r_{\omega\eta}}{r_\omega}=\frac{r_{\omega'\eta}}{r_\omega'}.
\end{equation*}
\end{itemize}
 If $\mathcal{G}$ fulfills the above Moran structure,
we call

\begin{equation*}
F=\bigcap\limits_{n=1}^\infty\bigcup\limits_{\omega\in D_n}B_\omega
\end{equation*}
the Moran set associated with $\mathcal{G}.$

For $k\in\mathbb{N},$ let

\begin{equation*}
c_k=\min\limits_{(i_1\cdots i_k)\in D_k}\frac{r_{i_1\cdots
i_k}}{r_{i_1\cdots i_{k-1}}},~~~M_k=\max\limits_{(i_1\cdots i_k)\in
D_k}r_{i_1\cdots i_k}
\end{equation*}

\begin{lem}{\rm\cite{FenLauWu}}\label{lem2.3}
For the Moran set $F$ defined as above, suppose furthermore

\begin{equation}\label{formular2.6}
\lim\limits_{k\to\infty}\frac{\log c_k}{\log M_k}=0.
\end{equation}
Then we have

\begin{equation*}
\dim_P F=\limsup\limits_{k\to\infty}s_k,
\end{equation*}
where $s_k$ satisfies the equation $\sum\limits_{\omega\in
D_k}r_\omega^{s_k}=1.$
\end{lem}

%%%%%%%%%%%%%%%%%%%%%%%%%%%%%%%%%%%%%%%%%%%%%%%%%%%%%%%%%%%%%%%%%%%%%%%%%%%%%%%%%%%%%%%%%%%%%%%%%%%%%%%%%%%%%%%%%%%
%%%%%%%%%%%%%%%%%%%%%%%%%%%%%%%%%%%%%%%%%%%%%%%%%%%%%%%%%%%%%%%%%%%%%%%%%%%%%%%%%%%%%%%%%%%%%%%%%%%%%%%%%%%%%%%%%%%
%%%%%%%%%%%%%%%%%%%%%%%%%%%%%%%%%%%%%%%%%%%%%%%%%%%%%%%%%%%%%%%%%%%%%%%%%%%%%%%%%%%%%%%%%%%%%%%%%%%%%%%%%%%%%%%%%%%
%%%%%%%%%%%%%%%%%%%%%%%%%%%%%%%%%%%%%%%%%%%%%%%%%%%%%%%%%%%%%%%%%%%%%%%%%%%%%%%%%%%%%%%%%%%%%%%%%%%%%%%%%%%%%%%%%%%
%%%%%%%%%%%%%%%%%%%%%%%%%%%%%%%%%%%%%%%%%%%%%%%%%%%%%%%%%%%%%%%%%%%%%%%%%%%%%%%%%%%%%%%%%%%%%%%%%%%%%%%%%%%%%%%%%%%

\section{Proof of Theorem \ref{thm1.5}}
To (1) in Theorem \ref{thm1.5}, the reader is referred to
\cite{LiWuXio} for a more detailed discussion. Let $K^I=\{x\in
K:A(D(x))\subseteq I\},$ and write $K_I=\{x\in K: A(D(x))=I\}.$ We
prove (2) by showing that

\begin{equation}\label{formula3.7}
\dim_PK_I\geq\sup\limits_{\alpha\in I}\beta^*(\alpha),
\end{equation}
and

\begin{equation}\label{formula3.8}
\dim_PK^I\leq\sup\limits_{\alpha\in I}\beta^*(\alpha).
\end{equation}
%%%%%%%%%%%%%%%%%%%%%%%%%%%%%%%%%%%%%%%%%%%%%%%%%%%%%%%%%%%%%%%%%%%%%%%%%%%%%%%%%%%%%%%%%%%%%%%%%%%%%%%%%%%%%%%%%%%
%%%%%%%%%%%%%%%%%%%%%%%%%%%%%%%%%%%%%%%%%%%%%%%%%%%%%%%%%%%%%%%%%%%%%%%%%%%%%%%%%%%%%%%%%%%%%%%%%%%%%%%%%%%%%%%%%%%
%%%%%%%%%%%%%%%%%%%%%%%%%%%%%%%%%%%%%%%%%%%%%%%%%%%%%%%%%%%%%%%%%%%%%%%%%%%%%%%%%%%%%%%%%%%%%%%%%%%%%%%%%%%%%%%%%%%
%%%%%%%%%%%%%%%%%%%%%%%%%%%%%%%%%%%%%%%%%%%%%%%%%%%%%%%%%%%%%%%%%%%%%%%%%%%%%%%%%%%%%%%%%%%%%%%%%%%%%%%%%%%%%%%%%%%
%%%%%%%%%%%%%%%%%%%%%%%%%%%%%%%%%%%%%%%%%%%%%%%%%%%%%%%%%%%%%%%%%%%%%%%%%%%%%%%%%%%%%%%%%%%%%%%%%%%%%%%%%%%%%%%%%%%

\subsection{Proof of inequality (\ref{formula3.7})}
Since $I$ is closed and $\beta^*(\alpha)$ is continuous, there
exists $\alpha_0\in I$ such that $\sup\limits_{\alpha\in
I}\beta^*(\alpha)=\beta^*(\alpha_0).$  The idea behind the proof is
to construct a Moran set $F$ such that

\begin{equation}\label{formula3.9}
F\subseteq K_I,
\end{equation}
and

\begin{equation}\label{formula3.10}
\dim_P F\geq\beta^*(\alpha_0).
\end{equation}
This approach was used in \cite{BaeOlsSni}, \cite{CawMau},
\cite{FenLau}, \cite{FenLauWu}, \cite{LiWuXio}, \cite{OlsWin} and
the following proof also benefits from these papers.

Let $i\in \mathbb{N}.$ Since $I$ is connected, we may choose
$q_{i,1}, \cdots, q_{i, M_i}\in \mathbb{R}$ such that

\begin{itemize}
  \item $\alpha_{i,j}\in I,$ where $\alpha_{i,j}=-\beta'(q_{i,j});$
  \item $I\subseteq \bigcup\limits_{j=1}^{M_i}B(\alpha_{i,j},\frac{1}{i});$
  \item $|\alpha_{i, j}-\alpha_{i, j+1}|\leq \frac{1}{i}$ for all
  $j,$~~~~$|\alpha_{i, M_i}-\alpha_{i+1,1}|\leq\frac{1}{i};$
  \item $\alpha_{i, M_i}=\alpha_0$ for all $i.$
\end{itemize}

\begin{rem}
In fact,
$\overline{\{\alpha_{i,j},\alpha_{i,j+1},\cdots,\alpha_{i,M_i},\alpha_{i+1,1},\alpha_{i+1,2}\cdots\}}=I,$
for any $i\in\mathbb{N}, 1\leq j\leq M_i.$
\end{rem}
Choose a positive sequence $(\delta_i)_{i=1}^\infty\searrow0.$ Note
that $\tau_{U_0}=\beta(q)$ for any $q\in \mathbb{R}.$ It follows
from Lemma \ref{lem2.2} that there exist positive real numbers
$(d_{i,j})_{i\in\mathbb{N}, j=1,\cdots, M_i},$

\noindent $(k_{i,j})_{i\in\mathbb{N}, j=1,\cdots, M_i}$ and
$\mathcal{B}_{i,j}=\{u_{i,j,s}:1\leq s\leq k_{i,j}\}\subset\Sigma^*$
such that\\

(i) $1>d_{1,1}>d_{1,2}>\cdots>d_{1,M_1}>d_{2,1}>d_{2,2}>\cdots>0.$\\

(ii)
$k_{i,j}\geq(d_{i,j})^{-\beta^*(\alpha_{i,j})+\delta_i(|q_{i,j}|+1)}.$\\

(iii) $(d_{i,j})^{1+\delta_i}\leq r_{u_{i,j,s}}\leq
(d_{i,j})^{1-\delta_i}$ for $1\leq s\leq k_{i,j}.$\\

(iv) $S_{u_{i,j,s}}(U(x_0,4r_0))(1\leq s\leq k_{i,j})$ are disjoint
subsets of $U_0.$\\

(v) $(d_{i,j})^{\alpha_{i,j}+3\delta_i}\leq p_{u_{i,j,s}}\leq
(d_{i,j})^{\alpha_{i,j}-3\delta_i}$ for $1\leq s\leq k_{i,j}.$\\

Choose a sequence of positive integers $(N_{k,j})_{k\in \mathbb{N},
j=1,\cdots, M_k}$ large enough such that

(vi) $(d_{k,j})^{N_{k,j}}\leq(d_{k, j+1})^{2^k}$ for any
$k\in\mathbb{N}$ and $1\leq j<M_k.$\\

(vii)
$\lim\limits_{k\to\infty}\frac{\sum\limits_{i=1}^{k-1}\sum\limits_{j=1}^{M_i}N_{i,j}\log
d_{i,j}+\sum\limits_{j=1}^{i_k-1}N_{k,j}\log d_{k,j}}{N_{k,i_k}\log
d_{k,i_k}}=0$ for any $1\leq i_k\leq M_k.$

Define a sequence of subsets of $\Sigma^*$ as follows:

\begin{equation*}
\underbrace{\mathcal{B}_{1,1},\cdots,\mathcal{B}_{1,1}}_{N_{1,1}{\rm
times}},
\underbrace{\mathcal{B}_{1,2},\cdots,\mathcal{B}_{1,2}}_{N_{1,2}{\rm
times}},\cdots
\underbrace{\mathcal{B}_{1,M_1},\cdots,\mathcal{B}_{1,M_1}}_{N_{1,M_1}{\rm
times}},
\underbrace{\mathcal{B}_{2,1},\cdots\mathcal{B}_{2,1}}_{N_{2,1}{\rm
times}},\cdots
\end{equation*}
and relabel them as $\left\{\mathcal{B}_n^*\right\}_{n=1}^\infty.$
Put

\begin{equation*}
\mathcal{G}=\left\{S_{v_1\cdots v_k}(\overline{U_0}):
k\in\mathbb{N},v_i\in\mathcal{B}_i^*{\rm ~~for~~}1\leq i\leq
k\right\},
\end{equation*}
and set

\begin{equation*}
F=\bigcap\limits_{n=1}^\infty\bigcup\limits_{v_1\in\mathcal{B}_1^*,
\cdots,v_n\in\mathcal{B}_n^*} S_{v_1\cdots v_n}(\overline{U_0}).
\end{equation*}
It is easy to check that $F$ is a Moran set associated with
$\mathcal{G},$ and $F$ satisfies (\ref{formula3.9}). The reader is
referred to \cite{LiWuXio} for a more detailed discussion.

For large $n,$ write
$n=\sum\limits_{j=1}^{M_1}N_{1,j}+\cdots+\sum\limits_{j=1}^{i_k}N_{k,j}+p$
with $1\leq p< N_{k,i_k+1}.$ Put $$A_k=\left(\prod\limits_{i=1}^{k-1}
\prod\limits_{j=1}^{M_i}(d_{i,j})^{N_{i,j}\left(-\beta^*(\alpha_{i,j})+\delta_i(|q_{i,j}|+1)\right)}\right)
\prod\limits_{j=1}^{i_k}(d_{k,j})^{N_{k,j}\left(-\beta^*(\alpha_{k,j})+\delta_k(|q_{k,j}|+1)\right)}.$$

Combining (ii) and (iii), we have

\begin{equation}
\prod\limits_{s=1}^n\sharp\mathcal{B}_s^*\geq
A_k(d_{k,i_k+1})^{p\left(-\beta^*(\alpha_{k,
i_k+1})+\delta_k(|q_{k,i_k+1}|+1)\right)}
\end{equation}
and

\begin{equation}\begin{split}\label{formula3.12}
&\inf\limits_{v_n\in\mathcal{B}_n^*}r_{v_n}\geq
(d_{k,i_k+1})^{1+\delta_k},\\&
\sup\limits_{v_1\in\mathcal{B}^*_1,\cdots,v_n\in\mathcal{B}_n^*}r_{v_1\cdots
v_n}\leq\left(\prod\limits_{i=1}^{k-1}\prod\limits_{j=1}^{M_i}(d_{i,j})^{N_{i,j}(1-\delta_i)}\right)\prod\limits_{j=1}^{i_k}
(d_{k,j})^{N_{k,j}(1-\delta_k)}(d_{k,i_k+1})^{p(1-\delta_k)}.
\end{split}\end{equation}
Using (\ref{formula3.12}) and (vi), we have

\begin{equation*}
\lim\limits_{n\to\infty}\frac{\log
\left(\inf\limits_{v_n\in\mathcal{B}_n^*}r_{v_n}\right)}{\log
\left(\sup\limits_{v_1\in\mathcal{B}_1^*,\cdots,v_n\in\mathcal{B}_n^*}r_{v_1\cdots
v_n}\right)}=0.
\end{equation*}
This implies that the condition (\ref{formular2.6}) in Lemma
\ref{lem2.3} is fulfilled. Hence, by Lemma \ref{lem2.3}, we conclude
that $\dim_P  F=\limsup\limits_{n\to\infty}s_n,$ where

\begin{equation*}
\sum\limits_{v_1\in\mathcal{B}_1^*,\cdots,v_n\in\mathcal{B}_n^*}(r_{v_1\cdots
v_n})^{s_n}=1.
\end{equation*}
It follows that

\begin{equation}\label{formula3.13}
\dim_P F\geq \limsup\limits_{n\to\infty}\frac{\log\left(
\prod\limits_{s=1}^n\sharp\mathcal{B}_s^*\right)}{-\log
(\inf\limits_{v_1\in\mathcal{B}_1^*,\cdots,v_n\in\mathcal{B}_n^*}r_{v_1\cdots
v_n})}.
\end{equation}
Next, let's prove (\ref{formula3.10}). Choose a special sequence of
positive integers $\{n_t\}_{t\geq1},$ with
$n_t=\sum\limits_{j=1}^{M_1}N_{1,j}+\cdots+\sum\limits_{j=1}^{M_t}N_{t,j}.$
Combining (ii) and (iii), we get
\begin{equation}\begin{split}\label{formula3.14}
\prod\limits_{s=1}^{n_t}\sharp\mathcal{B}^*_s&=k_{1,1}^{N_{1,1}}k_{1,2}^{N_{1,2}}
\cdots k_{1,M_1}^{N_{1,M_1}}k_{2,1}^{N_{2,1}}\cdots
k_{t,M_t}^{N_{t,M_t}}\\&\geq\prod\limits_{i=1}^t\prod\limits_{j=1}^{M_i}(d_{i,j})^{
N_{i,j}\left(-\beta^*(\alpha_{i,j})+\delta_i(|q_{i,j}|+1)\right)}.
\end{split}\end{equation}
This (\ref{formula3.14}), together with (\ref{formula3.13}) and the
following inequality

\begin{equation*}
\inf\limits_{v_1\in\mathcal{B}_1^*,\cdots,v_{n_t}\in\mathcal{B}_{n_t}^*}
r_{v_1\cdots
v_{n_t}}\geq\prod\limits_{i=1}^{t}\prod\limits_{j=1}^{M_i}(d_{i,j})^{N_{i,j}(1+\delta_i)},
\end{equation*}
yields

\begin{equation*}\begin{split}
\dim_PF&\geq\limsup\limits_{t\to\infty}
\frac{\log\prod\limits_{i=1}^t\prod\limits_{j=1}^{M_i}(d_{i,j})^{
N_{i,j}\left(-\beta^*(\alpha_{i,j})+\delta_i(|q_{i,j}|+1)\right)}}{-\log\prod\limits_{i=1}^{t}\prod\limits_{j=1}^{M_i}(d_{i,j})^{N_{i,j}(1+\delta_i)}}
\\&=\limsup\limits_{t\to\infty}
\frac{\sum\limits_{i=1}^t\sum\limits_{j=1}^{M_i}N_{i,j}
\left(-\beta^*(\alpha_{i,j})+\delta_i(|q_{i,j}|+1)\right)\log
d_{i,j}}{-\sum\limits_{i=1}^t\sum\limits_{j=1}^{M_i}N_{i,j}(1+\delta_i)\log
d_{i,j}}\\&\geq \beta^*(\alpha_0).
\end{split}\end{equation*}
%%%%%%%%%%%%%%%%%%%%%%%%%%%%%%%%%%%%%%%%%%%%%%%%%%%%%%%%%%%%%%%%%%%%%%%%%%%%%%%%%%%%%%%%%%%%%%%%%%%%%%%%%%%%%%%%%%%
%%%%%%%%%%%%%%%%%%%%%%%%%%%%%%%%%%%%%%%%%%%%%%%%%%%%%%%%%%%%%%%%%%%%%%%%%%%%%%%%%%%%%%%%%%%%%%%%%%%%%%%%%%%%%%%%%%%
%%%%%%%%%%%%%%%%%%%%%%%%%%%%%%%%%%%%%%%%%%%%%%%%%%%%%%%%%%%%%%%%%%%%%%%%%%%%%%%%%%%%%%%%%%%%%%%%%%%%%%%%%%%%%%%%%%%
%%%%%%%%%%%%%%%%%%%%%%%%%%%%%%%%%%%%%%%%%%%%%%%%%%%%%%%%%%%%%%%%%%%%%%%%%%%%%%%%%%%%%%%%%%%%%%%%%%%%%%%%%%%%%%%%%%%
%%%%%%%%%%%%%%%%%%%%%%%%%%%%%%%%%%%%%%%%%%%%%%%%%%%%%%%%%%%%%%%%%%%%%%%%%%%%%%%%%%%%%%%%%%%%%%%%%%%%%%%%%%%%%%%%%%%

\subsection{Proof of inequality (\ref{formula3.8})}
The idea of the following proof comes form an article by Patzschke
\cite{Pat}. First, present a lemma as below. Put $\Xi_{\leq
\alpha}=\left\{x\in
K:\limsup\limits_{r\to0}\frac{\log\mu(B(x,r))}{\log
r}\leq\alpha\right\}$ and write $\Xi_{\geq \alpha}=\left\{x\in
K:\liminf\limits_{r\to0}\frac{\log\mu(B(x,r))}{\log
r}\geq\alpha\right\}.$
\begin{lem}
For $\alpha\in [\alpha_{\min},\alpha_{\max}],$

\noindent{\rm(i)} if  $\alpha\leq\alpha(0),$ then $\dim_P\Xi_{\leq
\alpha}\leq\beta^*(\alpha).$

\noindent{\rm(ii)}if $\alpha\geq\alpha(0),$ then $\dim_P\Xi_{\geq
\alpha}\leq\beta^*(\alpha).$

\end{lem}

\noindent{\it Proof.} (i) If $\alpha\leq\alpha(0),$  then $q\geq0.$
Write $\beta=\beta(q).$  Let $\epsilon>0$ and $0<\rho<\frac{1}{2}$
and define

\begin{equation*}
\Xi_{<\alpha,m}=\left\{x\in
K:\rho^{n(\alpha+\epsilon)}\leq\mu(B(x,\rho^n)){\rm~~for~all~}n\geq
m\right\}
\end{equation*}
for $m\in\mathbb{N}.$ Then
$\Xi_{\leq\alpha}\subseteq\bigcup\limits_{m=1}^{\infty}\Xi_{<\alpha,m}.$

Let $B(x_1,r_1^*), B(x_2,r_2^*),\cdots$ be a $\rho^m$-packing of
$\Xi_{<\alpha,m}.$ For $i\in\mathbb{N}$ choose $n_i\in\mathbb{N}$
such that $\rho^{n_i}\leq r_i^*<\rho^{n_i}-1.$ Then $n_i\geq m$ and
$B(x_i,\rho^{n_i})\subseteq B(x_i,r_i^*)\subseteq
B(x_i,\rho^{n_i-1}).$ Hence, the sequence
$B(x_1,\rho^{n_1}),B(x_2,\rho^{n_2}),\cdots$ consists of disjoint
sets. For $0<r<1,$ define $\Gamma_r=\{u\in \Sigma^*: r_u<r\leq
r_{u^-}\}.$ It is well known that the OSC implies that
$\sup\limits_{x\in \mathbb{R}^d,0<r<1}\sharp\{u\in\Gamma_r:K_u\cap
B(x,r)\neq\emptyset\}<\infty.$ Then $\mu(B(x,r))\leq
c_1\max\{p_u:u\in\Gamma_r, K_u\cap B(x,r)\neq\emptyset\}$ for all
$n\in\mathbb{N}$ and all $x\in K.$ For $n\in \mathbb{N}$ write
$\Gamma(n)=\Gamma_{\rho^n}.$ By volume estimating we obtain a
constant $c_2$ such that

\begin{equation*}
\#\{i=1,2,\cdots:n_i=n, B(x_i,\rho^{n_i})\cap K_u\neq\emptyset\}\leq
c_2
\end{equation*}
for all $u\in\Gamma(n)$ and all $n.$ Therefore, using the definition
of $\Xi_{<\alpha, m},$

\begin{equation*}\begin{split}
&\sum\limits_{i=1}^\infty d(B(x_i, r^*_i))^{\alpha
q+\beta+\epsilon(1+q)}\\&\leq\left(\frac{2}{\rho}\right)^{\alpha
q+\beta+\epsilon(1+q)}\sum\limits_{i=1}^\infty\rho^{n_i(\beta+\epsilon)}\rho^{n_iq(\alpha+\epsilon)}
\\&\leq\left(\frac{2}{\rho}\right)^{\alpha
q+\beta+\epsilon(1+q)}\sum\limits_{i=1}^\infty\rho^{n_i(\beta+\epsilon)}\mu(B(x_i,\rho^{n_i}))^q
\\&\leq\left(\frac{2}{\rho}\right)^{\alpha
q+\beta+\epsilon(1+q)}c_1^qc_2\sum\limits_{n=m}^\infty\sum\limits_{u\in\Gamma(n)}\rho^{n(\beta+\epsilon)}p_u^q
\\&
\leq\left(\frac{2}{\rho}\right)^{\alpha
q+\beta+\epsilon(1+q)}c_1^qc_2
r_{\min}^{-1{(\beta+\epsilon)}}\sum\limits_{k=1}^{\infty}(\sum\limits_{k=1}^Np_i^qr_i^{\beta(q)+\epsilon})^k
\end{split}\end{equation*}
This shows, that $\dim_P\Xi_{<a,m}\leq \alpha
q+\beta(q)+\epsilon(1+q)$ for all $m\in\mathbb{N}$ and hence,
$\dim_P\Xi_{\leq \alpha}\leq\alpha q+\beta(q)+\epsilon(1+q).$ Since
$\epsilon>0$ was arbitrary, $\dim_P \Xi_{\leq\alpha}\leq\inf\{\alpha
q +\beta(q)\}.$

(ii) Let $\alpha\geq\alpha(0).$ Let $q<0$ and write
$\beta=\beta(q).$ Let $\epsilon>0$ and $0<\rho<\frac{1}{2}$ and
define

\begin{equation*}
\Xi_{>\alpha,m}=\left\{x\in
K:\rho^{n(\alpha-\epsilon)}\geq\mu(B(x,\rho^n)){\rm~~for~all~}n\geq
m\right\}
\end{equation*}
for $m\in\mathbb{N}.$ Then
$\Xi_{\geq\alpha}\subseteq\bigcup\limits_{m=1}^{\infty}\Xi_{>\alpha,m}.$
Let $B(x_1,r_1^*), B(x_2,r_2^*),\cdots$ be a $\rho^m$-packing of
$\Xi_{>\alpha,m}.$ For $i\in\mathbb{N}$ choose $n_i\in\mathbb{N}$
such that $\rho^{n_i}\leq r_i^*<\rho^{n_i}-1.$ Then $n_i\geq m$ and
$B(x_i,\rho^{n_i})\subseteq B(x_i,r_i^*)\subseteq
B(x_i,\rho^{n_i-1}).$ Hence, the sequence
$B(x_1,\rho^{n_1}),B(x_2,\rho^{n_2}),\cdots$ consists of disjoint
sets.   Using similar steps as above, then

\begin{equation*}\begin{split}
&\sum\limits_{i=1}^\infty d(B(x_i, r_i^*))^{\alpha
q+\beta+\epsilon(1-q)}\\&\leq\left(\frac{2}{\rho}\right)^{\alpha
q+\beta+\epsilon(1-q)}\sum\limits_{n=m}^\infty\sum\limits_{u\in\Gamma(n)}\rho^{n(\beta+\epsilon
)}p_u^q\\&\leq \left(\frac{2}{\rho}\right)^{\alpha
q+\beta+\epsilon(1-q)}r_{\min}^{-1(\beta+\epsilon)}\sum\limits_{k=1}^\infty(\sum\limits_{i=1}^N
p_i^qr_i^{\beta(q)+\epsilon})^k.
\end{split}\end{equation*}
Finally, it follows that
$\dim_P\Xi_{\geq\alpha}\leq\beta^*(\alpha).$ This completes the
proof of lemma.

If $\sup I\leq \alpha(0),$ then $K^I\subseteq\Xi_{\leq\sup_I}.$ So,
$\dim_P K^I\leq\sup\limits_{\alpha\in I}\beta^*(\alpha).$

If $\inf I\geq \alpha(0),$ then $K^I\subseteq\Xi_{\geq\inf_I}.$ So,
$\dim_P K^I\leq\sup\limits_{\alpha\in I}\beta^*(\alpha).$

If $\alpha(0)\in I,$ then $\dim_P K^I=\dim_H
K^I=\beta^*(\alpha(0)).$

%%%%%%%%%%%%%%%%%%%%%%%%%%%%%%%%%%%%%%%%%%%%%%%%%%%%%%%%%%%%%%%%%%%%%%%%%%%%%%%%%%%%%%%%%%%%%%%%%%%%%%%%%%%%%%%%%%%
%%%%%%%%%%%%%%%%%%%%%%%%%%%%%%%%%%%%%%%%%%%%%%%%%%%%%%%%%%%%%%%%%%%%%%%%%%%%%%%%%%%%%%%%%%%%%%%%%%%%%%%%%%%%%%%%%%%
%%%%%%%%%%%%%%%%%%%%%%%%%%%%%%%%%%%%%%%%%%%%%%%%%%%%%%%%%%%%%%%%%%%%%%%%%%%%%%%%%%%%%%%%%%%%%%%%%%%%%%%%%%%%%%%%%%%
%%%%%%%%%%%%%%%%%%%%%%%%%%%%%%%%%%%%%%%%%%%%%%%%%%%%%%%%%%%%%%%%%%%%%%%%%%%%%%%%%%%%%%%%%%%%%%%%%%%%%%%%%%%%%%%%%%%
%%%%%%%%%%%%%%%%%%%%%%%%%%%%%%%%%%%%%%%%%%%%%%%%%%%%%%%%%%%%%%%%%%%%%%%%%%%%%%%%%%%%%%%%%%%%%%%%%%%%%%%%%%%%%%%%%%%

\subsection{Proof of corollary \ref{corollary}}
 It follows from the following facts:

\begin{itemize}
  \item For any $\alpha\in I,$ we have that$ \left\{x\in K:A(D(x))=\alpha
\right\}\subseteq K^I. $

  \item $\dim_H K^I\leq\dim_PK^I.$
\end{itemize}

\begin{rem}
For $K^I,$ "I" can be any subset of $\mathbb{R}.$ In fact,
$K^I=K^{I\cap[\alpha_{\min},\alpha_{\max}]}.$
\end{rem}

%%%%%%%%%%%%%%%%%%%%%%%%%%%%%%%%%%%%%%%%%%%%%%%%%%%%%%%%%%%%%%%%%%%%%%%%%%%%%%%%%%%%%%%%%%%%%%%%%%%%%%%%%%%%%%%%%%%
%%%%%%%%%%%%%%%%%%%%%%%%%%%%%%%%%%%%%%%%%%%%%%%%%%%%%%%%%%%%%%%%%%%%%%%%%%%%%%%%%%%%%%%%%%%%%%%%%%%%%%%%%%%%%%%%%%%
%%%%%%%%%%%%%%%%%%%%%%%%%%%%%%%%%%%%%%%%%%%%%%%%%%%%%%%%%%%%%%%%%%%%%%%%%%%%%%%%%%%%%%%%%%%%%%%%%%%%%%%%%%%%%%%%%%%
%%%%%%%%%%%%%%%%%%%%%%%%%%%%%%%%%%%%%%%%%%%%%%%%%%%%%%%%%%%%%%%%%%%%%%%%%%%%%%%%%%%%%%%%%%%%%%%%%%%%%%%%%%%%%%%%%%%
%%%%%%%%%%%%%%%%%%%%%%%%%%%%%%%%%%%%%%%%%%%%%%%%%%%%%%%%%%%%%%%%%%%%%%%%%%%%%%%%%%%%%%%%%%%%%%%%%%%%%%%%%%%%%%%%%%%

\noindent {\bf Acknowledgements.}   The work was supported by the
National Natural Science Foundation of China (10971100) and National
Basic Research Program of China (973 Program) (2007CB814800).
%%%%%%%%%%%%%%%%%%%%%%%%%%%%%%%%%%%%%%%%%%%%%%%%%%%%%%%%%%%%%%%%%%%%%%%%%%%%%%%%%%%%%%%%%%%%%%%%%%%%%%%%%%%%%%%%%%%
%%%%%%%%%%%%%%%%%%%%%%%%%%%%%%%%%%%%%%%%%%%%%%%%%%%%%%%%%%%%%%%%%%%%%%%%%%%%%%%%%%%%%%%%%%%%%%%%%%%%%%%%%%%%%%%%%%%
%%%%%%%%%%%%%%%%%%%%%%%%%%%%%%%%%%%%%%%%%%%%%%%%%%%%%%%%%%%%%%%%%%%%%%%%%%%%%%%%%%%%%%%%%%%%%%%%%%%%%%%%%%%%%%%%%%%
%%%%%%%%%%%%%%%%%%%%%%%%%%%%%%%%%%%%%%%%%%%%%%%%%%%%%%%%%%%%%%%%%%%%%%%%%%%%%%%%%%%%%%%%%%%%%%%%%%%%%%%%%%%%%%%%%%%
%%%%%%%%%%%%%%%%%%%%%%%%%%%%%%%%%%%%%%%%%%%%%%%%%%%%%%%%%%%%%%%%%%%%%%%%%%%%%%%%%%%%%%%%%%%%%%%%%%%%%%%%%%%%%%%%%%%

\end{document}